\documentclass[12pt]{article}

\usepackage{amsmath}
\usepackage{amsfonts,graphicx}
\usepackage{amssymb,amsthm}
\newtheorem{theorem}{Theorem}
\newtheorem{corollary}{Corollary}
\newtheorem{lemma}{Lemma}
\newtheorem{definition}{Definition}
\newtheorem{proposition}{Proposition}
\newtheorem{example}{Example}
\numberwithin{equation}{section}

\begin{document}

\sloppy

\title{On the Stability of a Polling System with
an Adaptive Service Mechanism }

\author{Natalia Chernova\footnote{
\emph{Novosibirsk State University; 
E-mail: cher@nsu.ru}},
Sergey Foss\footnote{
\emph{ Heriot-Watt University,
Edinburgh and Institute of Mathematics, Novosibirsk; 
E-mail: s.foss@hw.ac.uk}},  
and Bara Kim\footnote{
\emph{
  Korea University, Seoul; 
E-mail: bara@korea.ac.kr}
}}

\date{}

\maketitle

\begin{abstract}
We consider a single-server cyclic polling system with three
queues where the server follows an adaptive rule: if it finds one
of queues empty in a given cycle, it decides not to visit that
queue in the next cycle. In the case of limited service policies,
we prove stability and instability results under some conditions
which are sufficient but not necessary, in general. Then we
discuss open problems with identifying the exact stability region
for models with limited service disciplines: we conjecture that a
necessary and sufficient condition for the stability may depend on
the whole distributions of the primitive sequences, and illustrate
that by examples. We conclude the paper with a section on the
stability analysis of a polling system with either gated or
exhaustive service disciplines.


\end{abstract}

\section{Introduction} \label{section1}

A standard polling model is a single-server system where the server
visits a finite  number of queues in cyclic order. The stability and
performance analysis of polling models with cyclic and other
scheduling policies has been a very popular research topic for
several decades. See, e.g. Borst (1995); Boxma et al. (2009);
Wierman et al. (2007); Winands et al. (2009) and the lists of
references therein for the progress in the studies of polling
models. One of the key tools in the modern stability analysis of
 queueing networks is the fluid limit approach. See, e.g.,
Chen and Mandelbaum (1991); Rybko and Stolyar (1992); Stolyar (1995);
Dai (1995); Dai and Meyn (1995); Down
(1996) for the detailed description. This approach involves the
use of the functional strong law of large numbers and works
perfectly well if one can identify all the limiting
parameters/functions there. This holds, in particular, for models
where the server's scheduling is state-independent.

In this paper, we make an attempt to study a model with an
adaptive scheduling of the server's visits to queues. Namely, we
consider a model with 3 queues and a cyclic policy, but assume in
addition that if in some cycle the server finds queue 2 empty, it
does not visit this queue in the next cycle. This may make sense
if one assumes, in addition, that the direct walking time from
queue 1 to queue 3 is smaller than that via queue 2. We consider
 three types of server's policies/disciplines in queues: limited,
 gated and exhaustive. We were motivated by the paper Vishnevsky and Semenova
(2008) where the authors considered a general model with many
queues and proposed a numerical iterative scheme to calculate the
mean waiting time. They showed numerically that, in the case of
limited service disciplines, the stability region may become
bigger by implementing such an adaptive scheduling.

Our original intention was to bring a mathematical accuracy here
and to obtain, for the proposed scheme,  some necessary and
sufficient conditions for the stability which are better than
those in non-adaptive schemes. We planned to use the fluid limit
approach for the stability. However, we have completed only a part
of our programme. First, we derived the evolution differential
equations for the fluid limits. Then we got two types of results.

When the service discipline in each queue is either gated or
exhaustive, the exact stability region is derived. The stability
region in this case is the same as for the standard non-adaptive
schemes, and the stability analysis here is a straightforward
application of the techniques developed, say, in Dai (1995) and
Dai and Meyn (1995). So, a possible advantage of the use of
adaptive scheduling is not in the stability, but in the
performance: when the polling systems are stable, stationary
characteristics under the adaptive
schemes may be smaller than those under the non-adaptive
schemes. 

If the server discipline in each queue is limited, then we can
derive tight fluid model equations only in the case where the
service discipline in the second queue is 1-limited. For a general
limited discipline, we use simple bounds for a fluid model
equation to obtain separately sufficient conditions for stability
and for instability. 

It is known that a sample-path analysis of fluid limits may become
inefficient due to several reasons. One of them was discussed, for 
example, by
Foss and Kovalevskii (1999) who considered fluid limits as weak
limits of stochastic processes under the linear scaling in time and in
space and show, in particular, that these limits may stay random.
We will discuss another reason for the limitation of the
fluid limit approach which is observed in the polling models:
a stability region may depend on
the entirely whole distribution of driving sequences.
We believe that the gap between the stability and instability
conditions is not just a technical limitation of the fluid
approximation, but indicates that the conditions obtained are
as sharp as one could get based only on the knowledge of the
first moments. We guess that, for any set of parameters in the gap,
one may propose a system that is stable and another system that is
unstable.

 The remainder of this paper is organized as follows.
 In Section \ref{section2}, we provide a detailed description of
the polling system and of the underlying Markov process. Then we
present the main results of the paper on stability and instability
of the system.  In Section \ref{section3}, we introduce fluid
limits for the system with limited service policies, derive
dynamical fluid equations and recall known stability and
instability criteria via fluid limits. Then we prove in Section
\ref{section4} the main results, in the case of limited service
policies. Section \ref{section5} presents simulation results and a
discussion on applicability and accuracy of the fluid model
approach for (in)stability of polling systems with adaptive
routing and limited service disciplines. In Section
\ref{section6}, we study a fluid model for gated and exhaustive
service disciplines and prove the main theorems in this case.

\section{Model Description and main results} \label{section2}

We consider a polling system with three infinite-buffer
queues/stations and a single server. The input stream to queue
$k=1,\,2,\,3$ is described by interarrival times $\tau_k(n)$ between
$n$th and $(n+1)$st customers, $n=1,\,2,\ldots$; and $n$th customer
in queue $k$ requires $\sigma_k(n)$ units of time for service.  The
server visits queues in a cyclic order. We assume that the cycles
start from queue 1, and that there are two possible cycle types:
 \begin{itemize}
 \item[$\bullet$] {\it standard} cycle: $1\to 2\to 3\to 1$;
 \item[$\bullet$] {\it reduced} cycle: $1\to 3\to 1$.
\end{itemize}
A choice of the type of the next cycle on the base of the current
cycle information is decribed in detail later on. In short, the
server may decide not to visit queue 2 in the next cycle if it finds
that queue empty.

The $n$th walking (switch-over) time from station $k$ to station
$(k+1) (\textrm{mod}\, 3)$ is denoted by $\xi_k(n)$, while walking
times from station 1 to station 3 are $\xi_4(n)$.

A service discipline describes the number of customers that (may be)
served during a single visit of the server to a queue. We will
assume that the service disciplines in all queues are the same and
belong to one of the following three classes:

\begin{enumerate}
 \item  {\it limited}\,: upon a visit to queue $k=1,\,2,\,3$, the server  serves at
most $l_k\ge 1$ customers, i.e. it continues working until either a
predefined number $l_k$ of customers is served or until the queue
becomes empty, whichever occurs first;
 \item  {\it gated}\,: upon
 a visit to a queue, the server serves all the customer that are
present at the moment of the arrival, and only those customers;
 \item {\it exhaustive}\,: the server continues to serve customers in a queue until
 the queue becomes empty.
\end{enumerate}

Now we describe inductively the way how does the server choose a
type of the next cycle. Here the state of queue 2 (empty or not)
within a cycle plays a crucial role.
 \begin{itemize}
 \item[$\bullet$] The first cycle is of the standard type.
 \item[$\bullet$] If cycle $n$ is of the standard type, and
 \begin{itemize}
 \item if, within this cycle, the server finds queue $2$ non-empty and proceeds with
 the service there, then cycle $n+1$ is again standard;
 \item if queue 2 is empty, the cycle $n+1$ is reduced.
 \end{itemize}
 \item[$\bullet$] If cycle $n$ is reduced, then cycle $n+1$ is standard.
 \end{itemize}

 \vspace{0.2cm}

 {\noindent \bf Stochastic Assumptions:}
We assume that all sequences of random variables introduced above
are mutually independent and that each of these sequences is i.i.d.,
with a finite first moment. We introduce arrival and service rates
at each queue as
\[
\lambda_k=\dfrac{1}{\mathsf E\tau_k(1)} >0, \qquad
\mu_k=\dfrac{1}{\mathsf E\sigma_k(1)} > 0.
\]
Let $\beta_k=\lambda_k/\mu_k$, and $\rho_0= \sum^3_{k=1} \beta_k$.
Further, we introduce switch-over rates as
\[
\nu_k=\dfrac{1}{\mathsf E \xi_k(1)}, \quad k=1,\ldots ,4.
\]
Then $\zeta=\sum_{k=1}^3\nu_k^{-1}$ is the mean switch-over time in
the standard cycle, and $\zeta^*=\nu_3^{-1}+\nu_4^{-1}$ is the mean
switch-over time in the reduced cycle.
 We assume that
\begin{equation}\label{eq21}
 \zeta > \zeta^*
\end{equation}
which makes sense: the server may decide to skip its visit to queue
2 if this may decrease the cycle time.

 \vspace{0.4cm}

 {\noindent \bf Representing Markov process:}
 We use
 \begin{eqnarray*}
 X(t)&=&(Q(t),\,A(t),\,B(t),\,B^0(t),\,H(t),\,I(t),\,C(t))
 \end{eqnarray*}
 to denote the state of our polling system at time $t$.
 The components are described below:
 \begin{itemize}
 \item[$\bullet$] $Q(t)=(Q_1(t),Q_2(t),Q_3(t))$, where $Q_i(t)$ is
 the total number of customers that are either waiting in queue $i$
 or being served at station $i$ at time $t$.
 \item[$\bullet$] $A(t)=(A_1(t),A_2(t),A_3(t))$, where $A_i(t)$ is
 the remaining interarrival time of the arrival process to
queue $i$ at time $t$.
 \item[$\bullet$] $B(t)$ is the remaining service time of a customer who is
 being served at time $t$ if the server is in service (and
 $B(t)=0$ if the server is walking at time $t$).
 \item[$\bullet$] $B^0(t)$ is the remaining walking time of the server
at time $t$,
 if the server is walking (and
 $B^0(t)=0$  if the server is in service).
 \item[$\bullet$] If the server is in service at time $t$, $H(t)$ is the station the server
 works at. If the server is walking at time $t$,
 $H(t)$ is the station the server is walking to.
 \item[$\bullet$] $I(t)$ takes values 0 and 1.
 It switches from 0 to 1 at the moment when the
 server arrives at queue 2 and finds it empty, and, vise versa, from
 1 to 0 when the server starts its walking from queue 1 to queue 3
in a reduced cycle.
 \item[$\bullet$] $C(t)$ is an additional parameter at time $t$ that makes
 $X(t)$ a Markov process.
If the server is in service at time $t$,
 then $C(t)$ is the number of service completions at queue $H(t)$ by time
$t$ during
 the current visit to the queue, for limited service disciplines; and
the
 number of arrivals at queue $H(t)$ by time $t$ during
 the current visit to the queue, for gated service disciplines.
In both cases we let $C(t)=0$ when the server is walking.
 In the case of exhaustive disciplines,
we do not need an extra parameter, so let $C(t)\equiv 0$.
 \end{itemize}
 The processes $X =\{X(t):t\ge 0\}$ are taken to be right-continuous
 with left limits.
 It follows from Dai (1995) that $X$ is a strong Markov process
 whose state space ${\cal S}$ is a subset of $\mathbb R_+^{11}$.

 The Markov process $X$ is said to be positive Harris recurrent if
 it possesses a unique stationary distribution.
 To state the main results of this paper, we make the following
 assumptions on the interarrival
 time distributions.
 For each $k=1,\,2,\,3$, we assume that the distribution of random variable
$\tau_k(1)$ has an unbounded support, i.e., ${\bf
P}(\tau_k(1)>t)>0$, for all $t>0$. Further, we assume that, for
each $k=1,\,2,\,3$, the distribution of $\tau_k(1)$ is spread-out,
i.e., there exists an integer $n>0$ and a non-negative function
$g(x)$ with $\int_0^\infty g(x)\,dx >0$, such that
$${\bf P}\left(a\leq\sum_{i=1}^n\tau_k(i)\leq b\right)\geq \int_a^b g(x)\,dx,
\;\text{ for any }\; 0\leq a < b.$$


We say that a polling system is {\it stable} if the underlying Markov process
is positive Harris recurrent, and {\it unstable}, otherwise. An unstable
polling model is {\it transient} if $|Q(t)|\to \infty$ a.s., as $t\to\infty$.

 The following theorems provides sufficient conditions for stability,
instability and transience of the
 polling systems under consideration.
\begin{theorem} \label{th1}
 A polling system with limited service disciplines
is stable if the following three
conditions hold:
\begin{eqnarray}
\rho_0+\dfrac{\lambda_2}{l_2}\zeta & < & 1, \label{st1}\\
\rho_0+\min\left(\dfrac{\lambda_1}{l_1}\dfrac{\zeta+\zeta^*}{2}+\lambda_2\dfrac{\zeta-\zeta^*}{2},\;\; \dfrac{\lambda_1}{l_1}\zeta \right) & < & 1, \label{st2}\\
\rho_0+\min\left(\dfrac{\lambda_3}{l_3}\dfrac{\zeta+\zeta^*}{2}+\lambda_2\dfrac{\zeta-\zeta^*}{2},\;\;
\dfrac{\lambda_3}{l_3}\zeta \right) & < & 1. \label{st3}
\end{eqnarray}
A polling system with either gated or exhaustive service
disciplines
is stable if 
$\rho_0 <1$.
\end{theorem}


\begin{theorem} \label{th1-1}
 A polling system with limited service disciplines
is unstable if at least one of the following three inequalities hold,
either
\begin{eqnarray}
\rho_0+\dfrac{\lambda_2}{l_2}\zeta & \ge & 1, \text{ or} \label{inst1-1} \\
\rho_0+\dfrac{\lambda_1}{l_1}\dfrac{\zeta+\zeta^*}{2}+\dfrac{\lambda_2}{l_2}\dfrac{\zeta-\zeta^*}{2} & \ge & 1, \text{ or} \label{inst2-1}\\
\rho_0+\dfrac{\lambda_3}{l_3}\dfrac{\zeta+\zeta^*}{2}+
\dfrac{\lambda_2}{l_2}\dfrac{\zeta-\zeta^*}{2}
& \ge & 1. \label{inst3-1}
\end{eqnarray}
Moreover, the system is transient if one of the inequalities above is strict.

A polling system with either gated or exhaustive service
disciplines 
is unstable if $\rho_0 \ge 1$, and transient if $\rho_{0} >1$.
\end{theorem}

In two particular cases, the statements of Theorems \ref{th1} and \ref{th1-1}
do match: 
\begin{corollary}\label{iff}
 A polling system with limited service disciplines and with $l_2=1$
 is stable if and only if
 (\ref{st1}), (\ref{st2}) and (\ref{st3}) are satisfied.\\
 A polling system with gated or exhaustive service
 disciplines is stable if 
 and only if $\rho_0<1$.
\end{corollary}


\section{The fluid model and (in)stability criteria -- Limited service disciplines} \label{section3}

 In this section, we consider a polling system with limited service
 disciplines.
 We define fluid limits and derive fluid model equations
 that are satisfied by fluid limits.
 Stability and instability criteria are given via the fluid
 model defined by the fluid model equations.
 In our construction of the fluid model, we follow the general
 scheme, see e.g. Chen and Mandelbaum (1991); Rybko and Stolyar (1992);
Stolyar (1995); Dai (1995); Dai and
 Meyn (1995); Down (1996).

 First we define processes related to our polling system.
 \begin{itemize}
 \item[$\bullet$] $T(t)=(T_1(t),T_2(t),T_3(t))$, $t \ge0$, where
 $T_k(t)$
 is the amount of time the server spends in service at station $k$
 during the time interval $[0,t]$.
 \item[$\bullet$] $U(t)=(U_1(t),U_2(t),U_3(t),U_4(t))$, $t \ge0$, where
 $U_1(t)$ ($U_2(t),U_3(t)$, and $U_4(t)$, respectively)
 is the amount of time the server spends walking from
 station 1 to 2 (2 to 3, 3 to 1, and 1 to 3, respectively),
during $[0,t]$.
 \end{itemize}
 Let ${\mathbb{X}}=\{(Q(t),T(t),U(t)):t\ge 0\}$.
 If the initial state $x\in {\cal S}$ of the Markov
 process $X$ is needed to be displayed explicitly,
 ${\mathbb{X}}_x$ is used for the process ${\mathbb{X}}$ obtained with the initial state $x$ of the Markov
 process $X$.

 By the strong law of large numbers, for almost all sample paths
 $\omega$, we have
 \begin{eqnarray}
 \lim_{n\to \infty} \frac{1}{n} \sum^n_{i=1} \tau_k(i,\omega) &=&
 \frac{1}\lambda_k, ~~~k=1,2,3, \label{slln-1} \\
 \lim_{n\to \infty} \frac{1}{n} \sum^n_{i=1} \sigma_k(i,\omega) &=&
 \frac{1}\mu_k,  ~~~k=1,2,3, \label{slln-2} \\
 \lim_{n\to \infty} \frac{1}{n} \sum^n_{i=1} \xi_k(i,\omega) &=&
 \frac{1}\nu_k,  ~~~k=1,2,3,4. \label{slln-3}
 \end{eqnarray}
It follows from the same argument as in Dai (1995) that for every
sample path $\omega$ satisfying (\ref{slln-1})-(\ref{slln-3}) and
every collection ${\{x_r:r>0\}}$ of initial states such that
 ${\{{|x_r|}/{r}:r>0\}}$ is bounded, there exists a
 subsequence $r_n\to \infty$ such that
 $\frac{1}{r_n}{{\mathbb{X}}}_{x_{r_n}}(r_n\cdot,\omega)$
 converges uniformly on any
 compact subset of $[0,\infty)$ to some limit say
 $ \bar{{\mathbb{X}}} =
 (\bar{Q}(\cdot),\bar{T}(\cdot),\bar{U}(\cdot))$.
 Each such limit $\bar{{\mathbb{X}}}$ is called a
 \emph{fluid limit}. In the special case where the sequence of initial
 states $\{x_r:r>0\}$ is independent of $r$, we call the limit a
 \emph{fluid limit with fixed initial state}. Both types of fluid limits
 are used in our subsequent stability analysis.

 It is known (see, e.g., Bramson (1998)) that, in the analysis of stability
via fluid limits, it
is sufficient to consider the so-called undelayed fluid limits only,
i.e. limits that satisfy the assumption
 \begin{eqnarray}
 \lim_{r \to \infty} \frac{1}{r} (|a_r|+b_r+b^0_r) &=& 0, \label{120}
 \end{eqnarray}
 where $a_r$, $b_r$ and $b^0_r$ are subvectors of the initial state $x_r=(q_r,a_r,b_r,b^0_r,h_r,i_r,c_r)$.
 From now on, we consider only undelayed fluid limits.

 If $\bar{\mathbb{X}}$ is a fluid limit obtained from
 a sequence of initial states $\{x_r\}$ satisfying (\ref{120}),
 then all of $\bar{Q}_k(\cdot)$, $k=1,2,3$,
 $\bar{T}_k(\cdot)$, $k=1,2,3$, and $\bar{U}_k(\cdot)$, $k=1,2,3,4$, are
Lipschitz continuous functions.
 Hence they are absolutely continuous and thus differentiable almost everywhere with
 respect to the Lebesgue measure.
 We say that $t$ is a {\it regular point} of
 $\bar{{\mathbb{X}}}$ if all components of
 $\bar{{\mathbb{X}}}$ are differentiable at $t$.
 In the rest of the paper, we implicitly assume that $t$
is a regular point whenever
 the derivative of a component of  $\bar{{\mathbb{X}}}$
 is involved.

 The following theorem presents equations that are satisfied by the fluid
 limits.
 \begin{lemma}\label{th2}
 For every fluid limit
 $\bar {\mathbb{X}}(t)=(\bar Q(t), \bar  T(t), \bar U(t))$,
 the following equations hold:
 \begin{eqnarray}
 && \bar Q_k(t)=\bar Q_k(0)+\lambda_k t - \mu_k \bar T_k(t), ~~
    k=1,2,3, ~ t \ge 0; \label{thm4-1} \\
 && \bar Q_k(t)\ge 0, ~~ k=1,2,3, ~ t \ge 0; \label{thm4-2}\\
 && \bar T_k(\cdot) \mbox{ and } \bar U_j(\cdot) \mbox{ are nondecreasing},
    ~~ k=1,2,3,~j=1,2,3,4;  \label{thm4-3}\\
 && \sum_{k=1}^3\bar T_k(t)+\sum_{j=1}^4\bar U_j(t)=t, ~~~ t \ge 0; \label{thm4-4} \\
 && \nu_1 \bar U_1(t) = \nu_2 \bar U_2(t), ~~ t \ge 0; \label{thm4-5}\\
 && \nu_3 \bar U_3(t) = \nu_1 \bar U_1(t) + \nu_{4} \bar U_{4}(t), ~~ t \ge 0;\label{thm4-6}\\
 && \nu_1 {\bar U}^\prime_1(t)+\nu_{4}{\bar U}^\prime_{4}(t)
    \ge  \dfrac{\mu_1}{l_1} {\bar T}^\prime_1(t), ~~ t \ge 0;  \label{thm4-7}\\
 && \nu_2 {\bar U}^\prime_2(t) \ge
    \dfrac{\mu_2}{l_2} {\bar T}^\prime_2(t), ~~ t \ge 0; \label{thm4-8}\\
 && \nu_3 {\bar U}^\prime_3(t) \ge
    \dfrac{\mu_3}{l_3} {\bar T}^\prime_3(t), ~~ t \ge 0; \label{thm4-9}\\
 && \frac{\mu_2}{l_2} \bar T_2'(t)\leq\nu_1{\bar U}^\prime_1(t)-\nu_4{\bar U}^\prime_4(t)\leq
    \mu_2 \bar T_2'(t), ~~ t \ge 0; \label{thm4-10}\\
 && \mbox{If } \bar Q_1(t)>0, \mbox{ then }\nu_1 {\bar U}^\prime_1(t)+\nu_{4}{\bar U}^\prime_{4}(t)
    = \dfrac{\mu_1}{l_1} {\bar T}^\prime_1(t),
    ~~ t\ge 0; \label{thm4-11}\\
 && \mbox{If } \bar Q_2(t)>0, \mbox{ then }  {\bar U}^\prime_4(t)=0
    \mbox{ and } \nu_2 {\bar U}^\prime_2(t) =
    \dfrac{\mu_2}{l_2} {\bar T}^\prime_2(t), ~~ t\ge 0;\label{thm4-12}\\
 && \mbox{If } \bar Q_3(t)>0, \mbox{ then } \nu_3 {\bar U}^\prime_3(t) =
    \dfrac{\mu_3}{l_3} {\bar T}^\prime_3(t), ~~ t\ge
    0.\label{thm4-13}
 \end{eqnarray}
 \end{lemma}
 \proof ~
 All equations in the theorem are obtained through the standard
 procedure.
 We provide below proofs only for equations (\ref{thm4-6}), (\ref{thm4-7}),
 (\ref{thm4-10}) and (\ref{thm4-11}), and omit all other
proofs.
 Let
 \begin{itemize}
 \item[$\bullet$] $D(t)=(D_1(t),D_2(t),D_3(t))$, where $D_i(t)$ is
 the number of service completions at station $i$ by time $t$;
 \item[$\bullet$] $M(t)=(M_1(t),M_2(t),M_3(t))$, where $M_i(t)$ is
 the number of service completions at station $i$
 if the server spends $t$ units of time working at station $i$;
 \item[$\bullet$] $E(t)=(E_1(t),E_2(t),E_3(t),E_4(t))$,
 where $E_1(t)$ ($E_2(t)$, $E_3(t)$, and $E_4(t)$, respectively) is
 the number of switch-over completions from station $1$ to station $2$
 (from $2$ to $3$,
 from $3$ to $1$, and
 from $1$ to $3$, respectively)
 by time $t$;
 \item[$\bullet$] $N(t)=(N_1(t),N_2(t),N_3(t),N_4(t))$,
 where $N_1(t)$ ($N_2(t)$, $N_3(t)$, and $N_4(t)$, respectively) is
 the number of switch-over completions from station $1$ to station $2$
 (from $2$ to $3$,
 from $3$ to $1$, and
 from $1$ to $3$, respectively)
 if the server spends $t$ units of time walking from station $1$ to station $2$
 (from $2$ to $3$,
 from $3$ to $1$, and
 from $1$ to $3$, respectively).
 \end{itemize}
 Then
 \begin{eqnarray}
 D_i(t)&=& M(T_i(t)), ~~ t\ge 0,~ i=1,2,3,\label{D}\\
 E_i(t)&=& N(U_i(t)), ~~ t\ge 0,~ i=1,2,3,4,\label{E}
 \end{eqnarray}
 and
 (\ref{slln-2}) and (\ref{slln-3}) imply
 \begin{eqnarray}
 \lim_{t \to \infty} \frac{M_i(t)}{t} &=& \mu_i,~ i=1,2,3, \label{lim-M} \\
 \lim_{t \to \infty} \frac{N_i(t)}{t} &=& \nu_i,~ i=1,2,3,4. \label{lim-N}
 \end{eqnarray}

 Now we prove (\ref{thm4-6}), (\ref{thm4-7}),
 (\ref{thm4-10}) and (\ref{thm4-11}).
 \begin{itemize}
 \item[$\bullet$] (\ref{thm4-6}):
 Since $|E_1(t)+E_4(t)-E_3(t)|\le 1$,
 (\ref{E}) yields
 $$|N_1(U_1(t))+N_4(U_4(t))-N_3(U_3(t))|\le 1.$$
 Applying fluid limits to the above equation, we obtain (\ref{thm4-6})
 with the help of (\ref{lim-N}).

 \item[$\bullet$] (\ref{thm4-7}):
 For $0\le t_1\le t_2$, we have
 $$E_1(t_2)-E_1(t_1)+E_4(t_2)-E_4(t_1)+1\ge
 \frac{D_1(t_2)-D_1(t_1)}{l_1}.$$
 Substituting
 (\ref{D}) and (\ref{E}) into the above equation yields
 $$N_1(U_1(t_2))-N_1(U_1(t_1))+N_4(U_4(t_2))-N_4(U_4(t_1))+1\ge \frac{M_1 (T_1(t_2))-M_1 (T_1(t_1))}{l_1}.$$
 Applying fluid limits to the above equation, we obtain (\ref{thm4-7})
 with the help of (\ref{lim-M}) and (\ref{lim-N}).

 \item[$\bullet$] (\ref{thm4-10}):
 Recall that if no customer is served at station 2 during a standard cycle,
then
 the next cycle is reduced. Conversely, if at least one customer is
 served at station 2 during a standard cycle, then the next cycle is
 standard.
 Therefore, for $0 \le t_1\le t_2$, the number of
 cycles with at least one service completion at station 2 during $(t_1,t_2)$
 differs at most by 1 from
 $(E_1(t_2)-E_1(t_1))-(E_4(t_2)-E_4(t_1))$.
 Therefore
 \begin{eqnarray*}
 \lefteqn{(E_1(t_2)-E_1(t_1))-(E_4(t_2)-E_4(t_1))-1} \\
 &\le& D_2(t_2)-D_1(t_1) \\
 &\le&  l_2((E_1(t_2)-E_1(t_1))-(E_4(t_2)-E_4(t_1))+1).
 \end{eqnarray*}
 Substituting
 (\ref{D}) and (\ref{E}) into the above equation and
 applying fluid limits yields
 (\ref{thm4-10}).

 \item[$\bullet$] (\ref{thm4-11}):
  Let $0\le t_1\le t_2 $.
 If $Q_1(t)>0$ for all $t\in(t_1,t_2)$, then
  $$\Big| E_1(t_2)-E_1(t_1)+E_4(t_2)-E_ 4(t_1)-
  \frac{D_1(t_2)-D_1(t_1)}{l_1}\Big|\le 1.$$
 Hence if $Q_1(t)>0$ for all $t\in(t_1,t_2)$, then
 $$\Big|N_1(U_1(t_2))-N_1(U_1(t_1))+N_4(U_4(t_2))-N_4(U_4(t_1))- \frac{M_1 (T_1(t_2))-M_1 (T_1(t_1))}{l_1}\Big|\le 1.$$
 Applying fluid limits to the above equation, we obtain
 (\ref{thm4-11}).

  \end{itemize}
  \qed
 \vspace{0.3cm}

 We call the equations (\ref{thm4-1})-(\ref{thm4-13}) the {\it fluid model
 equations} and call a solution $ \bar{\mathbb{X}} =
 \{(\bar{Q}(t),\bar{T}(t),\bar{U}(t)): t \ge 0\}$, of the
 fluid model equations a {\it fluid model solution}. Note that any fluid limit with
 fixed initial state necessarily has $\bar{Q}(0) = 0$. Thus these fluid limits form
 a subset of fluid model solutions with $\bar{Q}(0) = 0$. The following definitions and
 lemmas indicate the usefulness of different types of fluid limits.

 \begin{definition} \label{def1}
 \begin{enumerate}
 \item[(i)] The fluid model is {\it stable} if there exists a $\delta>0$ such that
 for each fluid model solution $ \bar{\mathbb{X}}$ with
 $|\bar{Q}(0)| \le 1$, $\bar Q(t)=0$ for $t \ge \delta$.
 \item[(ii)] The fluid model is {\it weakly unstable} if there exists a $\delta>0$
 such that, for each fluid model solution $ \bar{\mathbb{X}}$ with
 $\bar{Q}(0) = 0$, $\bar Q(\delta) \neq 0$.
 \end{enumerate}
 \end{definition}
The reasoning used in Dai (1995,1996), can be applied easily to our
polling system to obtain the following criteria.
 \begin{lemma} \label{lem-a1}
 (Dai (1995)) If the fluid model is stable, then the stochastic
polling system is stable too.
 \end{lemma}
 \begin{lemma} \label{lem-a2}
 (Dai (1996)) If the fluid model is weakly unstable, then
the stochastic polling system is transient.
 \end{lemma}

We present now a weaker instability criterion than Lemma \ref{lem-a2},
by
 applying to our polling systems
 the arguments first introduced in Dai et al. (2007), see Lemma
 \ref{lem-a3}
 below.
 If we assume \emph{a priori} that the process $X$ is positive
 Harris
 recurrent, then any fluid limit with fixed initial state must obey
 an extra dynamical equation, which augments the fluid model
 equations presented in (\ref{thm4-1})-(\ref{thm4-13}).
 Let $F_i(t)$ be the number of server's visits to station $i$ when
 it is empty, in the time interval $(0,t)$.
 By the theory of Markov regenerative processes, if $X$ is positive Harris
 recurrent, then there are positive numbers $f_i$, $i=1,2,3$, such
 that
 \begin{eqnarray}
 \lim_{t\to\infty} \frac{F_i(t)}{t} &=& f_i, ~~ i=1,2,3,\label{f_i}
 \end{eqnarray}
 with probability 1.

  \begin{lemma}\label{th2-aug}
  Suppose that the Markov process $X$ is positive
 Harris
 recurrent, and
 a fluid limit with fixed initial state $\bar {\mathbb{X}}(t)=(\bar Q(t), \bar  T(t), \bar U(t))$
 is driven from a sample path $\omega$ that satisfies (\ref{f_i}).
 Then the following inequalities hold: For $t\ge 0$,
 \begin{eqnarray}
  \nu_1 {\bar U}^\prime_1(t)+\nu_{4}{\bar U}^\prime_{4}(t)
    &>&  \dfrac{\mu_1}{l_1} {\bar T}^\prime_1(t); \label{aug1}\\
  \nu_2 {\bar U}^\prime_2(t) &>&
    \dfrac{\mu_2}{l_2} {\bar T}^\prime_2(t);\label{aug2}\\
 \nu_3 {\bar U}^\prime_3(t) &>&
    \dfrac{\mu_3}{l_3} {\bar T}^\prime_3(t).\label{aug3}
 \end{eqnarray}
 \end{lemma}
 \proof ~
 We have
 \begin{eqnarray*}
 D_1(t_2)-D_1(t_1)
 &\le & l_1 ((E_1(t_2)-E_1(t_1))+(E_4(t_2)-E_4(t_1))-(F_1(t_2)-F_1(t_1))+1), \\
 D_2(t_2)-D_2(t_1)
 &\le & l_2 ((E_2(t_2)-E_2(t_1))-(F_2(t_2)-F_2(t_1))+1), \\
 D_3(t_2)-D_3(t_1)
 &\le & l_3 ((E_3(t_2)-E_3(t_1))-(F_3(t_2)-F_3(t_1))+1).
 \end{eqnarray*}
 Substituting (\ref{D}) and (\ref{E}) into the above equations
 and applying fluid limits leads to
 \begin{eqnarray*}
 \mu_1 \bar T'_1(t)
 &\le & l_1 (\nu_1 \bar U'_1(t)+\nu_4 \bar U'_4(t)-f_1), \\
  \mu_2 \bar T'_2(t)
 &\le & l_2 (\nu_2 \bar U'_2(t)-f_2), \\
   \mu_3 \bar T'_3(t)
 &\le & l_3 (\nu_3 \bar U'_3(t)-f_3).
 \end{eqnarray*}
 Since $f_i>0$, $i=1,2,3$, Lemma \ref{th2-aug} is proved.
 \qed

 \vspace{0.3cm}

  We call the union of two systems of equations and inequalities
  (\ref{thm4-1})-(\ref{thm4-13})
  and (\ref{aug1})-(\ref{aug3})
  the {\it augmented fluid model
 equations} and call a solution $ \bar{\mathbb{X}}$, to these union
 an {\it augmented fluid model solution}.

 \begin{definition} \label{def2}
 The augmented fluid model is {\it weakly unstable} if there exists a $\delta>0$
 such that for each augmented fluid model solution $ \bar{\mathbb{X}}$, with
 $\bar{Q}(0) = 0$, $\bar Q(\delta) \neq 0$.
 \end{definition}

 Suppose that the augmented fluid model is weakly unstable but the
 Markov process $X$ is positive Harris recurrent.
 Since the augmented fluid model equations
 are satisfied by every fluid limit which is
 a limit of scaled sample paths with fixed
 initial state, the argument in Dai (1996) implies that
 the process is
 transient in the sense that, $|Q(t)|\to \infty$ as $ t \to \infty$ with
 probability 1, which is a contradiction.
 Therefore we obtain the following instability criterion.

 \begin{lemma} \label{lem-a3}
 (Dai et al. (2007))
 If the augmented fluid model is weakly unstable, then the stochastic
system is unstable.
 \end{lemma}

\section{Proof of Theorems \ref{th1} and \ref{th1-1} for limited service disciplines} \label{section4}

 For a polling system with  limited service
 disciplines, Theorems \ref{th1} and \ref{th1-1} 
follow from Propositions \ref{prop1}, \ref{prop2} and \ref{prop3}
 below and the (in)stability criteria, Lemmas \ref{lem-a1}, \ref{lem-a2} and \ref{lem-a3}.

 \begin{proposition} \label{prop1}
 For the polling system with limited service
 disciplines, the fluid model is stable if (\ref{st1}), (\ref{st2}) and
 (\ref{st3}) are satisfied.
 \end{proposition}

 \begin{proposition} \label{prop2}
 For the polling system with limited service
 disciplines, the augmented fluid model is weakly unstable if (\ref{inst1-1}), (\ref{inst2-1})
 or (\ref{inst3-1}) holds.
 \end{proposition}

 \begin{proposition} \label{prop3}
 For the polling system with limited service
 disciplines, the fluid model is weakly unstable if at least one of the
inequalities
(\ref{inst1-1}), (\ref{inst2-1})
 or (\ref{inst3-1}) is strict.
 \end{proposition}

 The remainder of this section is devoted to the proof of
 Propositions \ref{prop1}, \ref{prop2} and \ref{prop3}.
 For a fluid model solution $\bar{\mathbb{X}}$,
 let
 \begin{eqnarray}
 J(t)&\equiv & \{k\in \{1,2,3\}:\bar Q_k(t)>0\}, ~~ t\ge 0.
 \label{J(t)}
 \end{eqnarray}

 \begin{lemma}\label{lem4}
 For each fluid model solution $\bar{\mathbb{X}}$,
 if $\bar Q_2(t)>0$, then
 \begin{eqnarray}
 {\bar T}^\prime_j(t) &=& \dfrac{l_j}{\mu_j}
 \dfrac{1-\rho_0+\sum_{k\in J(t)}\beta_k}{\zeta+\sum_{k\in J(t)}
 \frac{l_k}{\mu_k}},~~ j \in J(t). \label{lem4-1}
 \end{eqnarray}
 \end{lemma}
 \proof
 Let $\bar{\mathbb{X}}$ be a fluid model solution. Suppose that $\bar Q_2(t)>0$.
 If $\bar Q_j(t)>0$, then,
 according to (\ref{thm4-1}), (\ref{thm4-5}), (\ref{thm4-6}) and (\ref{thm4-11})--(\ref{thm4-13}), we
 have
 \begin{eqnarray*}
 \bar U'_k(t) &=& \frac{\mu_j}{l_j} \frac{1}{\nu_k} \bar T'_j(t),
 ~~ k=1,2,3,\\
 \bar T'_k(t) &=& \left\{
 \begin{array}{cl}
 \frac{l_k}{\mu_k} \frac{\mu_j}{l_j}   \bar T'_j(t), & k \in J(t),
 \\
 \beta_k, & k \in \{1,2,3\}\setminus J(t).
 \end{array}
 \right.
 \end{eqnarray*}
 Substituting the above equations into (\ref{thm4-4}) yields
 (\ref{lem4-1}).
 \qed

 \begin{lemma}\label{lem5}
 For each fluid model solution $\bar{\mathbb{X}}$,
 if $\bar Q_2(t)=0$, then
 \begin{eqnarray}
 {\bar T}^\prime_j(t)
 &\ge & \dfrac{l_j}{\mu_j}
 \dfrac{1-\rho_0+\sum_{k\in J(t)}\beta_k-\lambda_2\frac{\zeta-\zeta^*}{2}}
 {\frac{\zeta+\zeta^*}{2}+\sum_{k\in J(t)}\frac{l_k}{\mu_k}},~~ j \in J(t); \label{lem5-1-1} \\
 {\bar T}^\prime_j(t)
 &\ge & \dfrac{l_j}{\mu_j}
 \dfrac{1-\rho_0+\sum_{k\in J(t)}\beta_k} {\zeta+\sum_{k\in J(t)}\frac{l_k}{\mu_k}},~~ j \in J(t). \label{lem5-2}
 \end{eqnarray}
 \end{lemma}
  \proof
 Let $\bar{\mathbb{X}}$ be a fluid model solution. Suppose that $\bar Q_2(t)=0$.
According to (\ref{thm4-1}), (\ref{thm4-5}), (\ref{thm4-6}) and (\ref{thm4-11})--(\ref{thm4-13}), we
have
 \begin{eqnarray}
 \bar T'_k(t) &=& \frac{l_k}{\mu_k} \frac{\mu_j}{l_j} \bar T'_j(t),
 ~~ \mbox{ if } \bar Q_j(t)>0 \mbox{ and } k \in J(t), \label{thm5-1}  \\
 \bar T'_k(t) &=& \beta_k, ~~ \mbox{ if } k \in \{1,2,3\}\setminus
 J(t).  \label{thm5-3}
 \end{eqnarray}
Clearly, for any fixed $k$, the right-hand side in \eqref{thm5-1} is 
the same for all $j$ with $\bar Q_j(t)>0$.

According to (\ref{thm4-5})-(\ref{thm4-9}) and (\ref{thm4-11})-(\ref{thm4-13}), we
have
 \begin{eqnarray}
 \sum^4_{k=1} \bar U'_k(t) &=& \frac{\mu_j}{l_j} \bar T'_j(t) \zeta-\nu_4 \bar U'_4(t) (\zeta-\zeta^*),
 ~~ \mbox{ if } \bar Q_j(t)>0.   \label{thm5-4}
 \end{eqnarray}
 By (\ref{thm4-3}), (\ref{thm4-6})-(\ref{thm4-11}) and
 (\ref{thm5-4}),
 \begin{eqnarray}
 \sum^4_{k=1} \bar U'_k(t) &\le& \frac{\mu_j}{l_j} \bar T'_j(t) \frac{\zeta+\zeta^*}{2}+
  \lambda_2 \frac{\zeta-\zeta^*}{2},
 ~~ \mbox{ if } \bar Q_j(t)>0,   \label{thm5-7}  \\
 \sum^4_{k=1} \bar U'_k(t) &\le& \frac{\mu_j}{l_j} \bar T'_j(t) \zeta,
 ~~ \mbox{ if } \bar Q_j(t)>0.   \label{thm5-8}
 \end{eqnarray}
 Substituting (\ref{thm5-1}), (\ref{thm5-3}) and (\ref{thm5-7}) into
 (\ref{thm4-4}) yields (\ref{lem5-1-1}).
 Substituting (\ref{thm5-1}), (\ref{thm5-3}) and (\ref{thm5-8}) into
 (\ref{thm4-4}) yields (\ref{lem5-2}).
 \qed

 \begin{lemma}\label{lem5-1}
 For each fluid model solution $\bar{\mathbb{X}}$,
 \begin{eqnarray}
 {\bar T}^\prime_2(t) &\le& \dfrac{l_2}{\mu_2}
 \dfrac{1-\rho_0+\sum_{k\in J(t)\cup\{2\}}\beta_k}{\zeta+\sum_{k\in J(t)\cup\{2\}}
 \frac{l_k}{\mu_k}};  \label{lem5-4}\\
 {\bar T}^\prime_j(t) &\leq& \dfrac{l_j}{\mu_j}
 \dfrac{1-\rho_0+\sum_{k\in J(t)\cup\{j\}}\beta_k-\frac{\lambda_2}{l_2}\frac{\zeta-\zeta^*}{2}}
 {\frac{\zeta+\zeta^*}{2}+\sum_{k\in J(t)\cup\{j\}}\frac{l_k}{\mu_k}},~~ j =1,3. \label{lem5-5}
 \end{eqnarray}
 \end{lemma}
  \proof
 Let $\bar{\mathbb{X}}$ be a fluid model solution.
 According to (\ref{thm4-1}), (\ref{thm4-5}), (\ref{thm4-6}) and (\ref{thm4-11})-(\ref{thm4-13}), we
 have
 \begin{eqnarray}
 \bar T'_k(t) &\ge& \frac{l_k}{\mu_k} \frac{\mu_j}{l_j} \bar T'_j(t),
 ~~   k \in J(t), ~j=1,2,3, \label{thm5-1-1}   \\
 \bar T'_k(t) &=& \beta_k, ~~   k \in \{1,2,3\}\setminus
 J(t).  \label{thm5-3-1}
 \end{eqnarray}
 According to (\ref{thm4-5})-(\ref{thm4-9}) and (\ref{thm4-11})-(\ref{thm4-13}), we
 have
 \begin{eqnarray}
 \sum^4_{k=1} \bar U'_k(t) &\ge& \nu_3 \bar U'_3(t)
 \zeta-\nu_4 \bar U'_4(t) (\zeta-\zeta^*),
 ~~ j=1,2,3,   \label{thm5-5}  \\
 \sum^4_{k=1} \bar U'_k(t) &\ge& \frac{\mu_2}{l_2} \bar T'_2(t) \zeta.  \label{thm5-6}
 \end{eqnarray}
 By (\ref{thm4-3}), (\ref{thm4-6})-(\ref{thm4-11}), (\ref{thm4-13})
 and (\ref{thm5-5}),
 \begin{eqnarray}
 \sum^4_{k=1} \bar U'_k(t) &\ge& \frac{\mu_j}{l_j} \bar T'_j(t)
\frac{\zeta+\zeta^*}{2}+
 \frac{\lambda_2}{l_2}  \frac{\zeta-\zeta^*}{2}.  \label{thm5-9}
 \end{eqnarray}
 Substituting (\ref{thm5-1-1}), (\ref{thm5-3-1}) and (\ref{thm5-6}) into
 (\ref{thm4-4}) yields (\ref{lem5-4}).
 Substituting (\ref{thm5-1-1}), (\ref{thm5-3-1}) and (\ref{thm5-9}) into
 (\ref{thm4-4}) yields (\ref{lem5-5}).
 \qed

 \begin{lemma}\label{lem6}
 For each augmented fluid model solution $\bar{\mathbb{X}}$,
 \begin{eqnarray}
 {\bar T}^\prime_2(t) &<& \dfrac{l_2}{\mu_2}
 \dfrac{1-\rho_0+ \beta_2}{\zeta+\frac{l_2}{\mu_2}};
 \label{lem6-1} \\
 {\bar T}^\prime_j(t) &<& \dfrac{l_j}{\mu_j}
 \dfrac{1-\rho_0+ \beta_j-\frac{\lambda_2}{l_2}\frac{\zeta-\zeta^*}{2}}
 {\frac{\zeta+\zeta^*}{2}+ \frac{l_j}{\mu_j}},~~ j=1,3. \label{lem6-2}
 \end{eqnarray}
 \end{lemma}
 \proof
 Let $\bar{\mathbb{X}}$ be an augmented fluid model solution.
 By (\ref{thm4-11})-(\ref{thm4-13}) and (\ref{aug1})-(\ref{aug3}),
 we have $\bar Q_k(t)=0$, $k=1,2,3$, and, by (\ref{thm4-1}),
 \begin{eqnarray}
 \bar T'_k(t)&=&\beta_k, ~~t\ge 0, ~ k=1,2,3.
 \label{lem6-3}
 \end{eqnarray}
 By (\ref{thm4-5}) and (\ref{thm4-6}) and (\ref{aug1})-(\ref{aug3}),
 \begin{eqnarray}
 \sum^4_{k=1} \bar U'_k(t)&>&\frac{\mu_2}{l_2} \bar T'_2(t) \zeta,
 \label{lem6-4} \\
 \sum^4_{k=1} \bar U'_k(t) &>& \frac{\mu_j}{l_j} \bar T'_j(t) \zeta-\nu_4 \bar U'_4(t) (\zeta-\zeta^*),
 ~~ j=1,3.
 \label{lem6-5}
 \end{eqnarray}
 Substituting (\ref{lem6-3}) and (\ref{lem6-4}) into (\ref{thm4-4})
 leads to (\ref{lem6-1}).
  Substituting (\ref{lem6-3}) and (\ref{lem6-5}) into (\ref{thm4-4})
 leads to (\ref{lem6-2}).
 \qed

 \vspace{0.2cm}

 Now we are ready to prove Propositions \ref{prop1}, \ref{prop2} and \ref{prop3}.

 \vspace{0.2cm}
 \noindent {\it Proof of Proposition \ref{prop1}.}
 ~ Suppose that (\ref{st1}), (\ref{st2}) and
 (\ref{st3}) are satisfied.
 For a fluid model solution
 $\bar{\mathbb{X}}$, let
 \begin{eqnarray*}
 W(t)&\equiv&\sum^3_{k=1}\frac{\bar Q_k(t)}{\mu_k}.
 \end{eqnarray*}
 Then $W(t)=0$ if and only if $|\bar Q(t)|=0$.
 We prove that the fluid model is stable
 by showing that there is a $\delta>0$ such that
 for each fluid model solution $ \bar{\mathbb{X}}$, with
 $|\bar{Q}(0)| \le 1$, $W(t)=0$ for $t \ge \delta$.
 The proof proceeds through 3 steps.

 \vspace{0.2cm}
 \noindent{\it Step 1.
 Let $m$ be the number of indices $j$ in $\{1,2,3\}$ such that
 $\frac{\lambda_j}{l_j}\le\frac{\lambda_2}{l_2}$, and let
 $j_1,\ldots,j_m$ be the indices enumerated so that
 $\frac{\lambda_{j_1}}{l_{j_1}}\le  \frac{\lambda_{j_2}}{l_{j_2}}\le  \ldots \le \frac{\lambda_{j_m}}{l_{j_m}}$.
 Then there exist $\delta_k\ge 0$, $1\le k \le m$, such that for each fluid model solution $\bar {\mathbb{X}}$
 with $|\bar Q(0)|\le 1$,
 \begin{eqnarray}
 \bar Q_j(t)=0 \mbox{ for } t \ge \delta_k \mbox{ and }j \in \{j_1,\ldots,j_k\}.
 \label{step1-1}
 \end{eqnarray}}
 \noindent{\it Proof. }
 Let $\delta_0=0$.
 For $k=0$, (\ref{step1-1}) holds trivially for each fluid model
 solution $\bar {\mathbb{X}}$ with $|\bar Q(0)|\le 1$.
 For $1\le k \le m$, suppose that there exists $\delta_{k-1}\ge 0$
 such that, for each fluid model solution $\bar {\mathbb{X}}$
 with $|\bar Q(0)|\le 1$,
 $\bar Q_j(t)=0$ for $t \ge \delta_{k-1}$ and $j \in \{j_1,\ldots,j_{k-1}\}$.
 Suppose that $ \bar{\mathbb{X}}$ is a fluid model solution with
 $|\bar{Q}(0)| \le 1$.
 According to (\ref{lem4-1}) and (\ref{lem5-2}), if $\bar
 Q_{j_k}(t)>0$, then
 \begin{eqnarray}
 \bar Q_{j_k}'(t)
 &\le& \lambda_{j_k}-l_{j_k} \frac{1-\rho_0+\sum_{i\in J(t)}\beta_i}{
       \zeta+\sum_{i\in J(t)}\frac{l_i}{\mu_i}} \nonumber \\
 &=& \frac{\lambda_{j_k} \zeta+\lambda_{j_k}\sum_{i\in J(t)}\frac{l_i}{\mu_i}-l_{j_k}
     (1-\rho_0+\sum_{i\in J(t)}\beta_i)}{\zeta+\sum_{i\in
     J(t)}\frac{l_i}{\mu_i}}. \label{step1-2}
 \end{eqnarray}
 Since $\{j_1,\ldots,j_{k-1}\}\cap J(t)=\phi$ for $t\ge
 \delta_{k-1}$,
 (\ref{step1-2}) leads to
 \begin{eqnarray*}
 Q_{j_k}'(t)
 &\le&  -\frac{l_{j_k}(1-\rho_0-\frac{\lambda_{j_k}}{l_{j_k}}\zeta)}{\zeta+\sum_{i\in
     J(t)}\frac{l_i}{\mu_i}}, ~~ \mbox{ if } t\ge \delta_{k-1}
     \mbox{ and } \bar Q_{j_k}(t)>0.
 \end{eqnarray*}
 Hence
 \begin{eqnarray*}
 \bar Q_{j_k}'(t)
 &\le&  -\epsilon_k, ~~ \mbox{ if } t\ge \delta_{k-1}
     \mbox{ and } \bar Q_{j_k}(t)>0,
 \end{eqnarray*}
 where
 \begin{eqnarray*}
 \epsilon_k
 &\equiv&  \frac{l_{j_k}(1-\rho_0-\frac{\lambda_{j_k}}{l_{j_k}}\zeta)}{\zeta+\sum^3_{i=1}\frac{l_i}{\mu_i}} >0.
 \end{eqnarray*}
 According to (\ref{thm4-1}) and  (\ref{thm4-3}),
 we have
  \begin{eqnarray*}
 Q_{j_k}(\delta_{k-1})
 &\le& 1+\lambda_{j_k}\delta_{k-1}.
 \end{eqnarray*}
 Hence $\bar Q_{j_k}(t)=0$ for $t \ge \delta_k$, where
 $\delta_k\equiv \delta_{k-1}+\frac{1+\lambda_{j_k}\delta_{k-1}}{\epsilon_k}$.
 The proof is completed by induction on $k$.
 \qed

 \vspace{0.2cm}
 \noindent{\it Step 2.  There exists an $\epsilon>0$ such that, for each fluid model solution $\bar
 {\mathbb{X}}$,
 $W'(t)\le -\epsilon$ if $W(t)>0$ and $\bar Q_2(t)=0$.}

 \vspace{0.2cm}
 \noindent{\it Proof. } It is easily proved that at least one of the
 following inequalities holds:
 \begin{eqnarray}
\rho_0+\dfrac{\lambda_j}{l_j}\zeta & < & 1, ~~ j=1,3,\label{step2-1}\\
\rho_0+\dfrac{\lambda_j}{l_j}\dfrac{\zeta+\zeta^*}{2}+\lambda_2\dfrac{\zeta-\zeta^*}{2}
& < & 1,~~ j=1,3.\label{step2-2}
\end{eqnarray}
 Let  $ \bar{\mathbb{X}}$ be a fluid model solution. Suppose that  $W(t)>0$ and $\bar
 Q_2(t)=0$.
 By (\ref{thm4-1}),
 \begin{eqnarray}
 W'(t)&=& \sum_{k \in J(t)}\beta_k -\sum_{k \in J(t)} \bar T'_k(t).
 \label{step2-3}
 \end{eqnarray}

 First suppose that (\ref{step2-1}) holds.
 Substituting (\ref{lem5-2}) into (\ref{step2-3}) leads to
 \begin{eqnarray*}
 W'(t)&\le& \frac{\sum_{k \in J(t)}\frac{l_k}{\mu_k}(\rho_0+\frac{\lambda_k}{l_k}\zeta-1)}{\zeta+\sum_{k\in J(t)}\frac{l_k}{\mu_k}}.
 \end{eqnarray*}
 Hence $W'(t)\le -\epsilon$, where
 $$\epsilon \equiv \min_{K\in \{\{1\},\{3\},\{1,3\}\}}
 \frac{\sum_{k \in K}\frac{l_k}{\mu_k}(1-\rho_0-\frac{\lambda_k}{l_k}\zeta)}{\zeta+\sum_{k\in
 K}\frac{l_k}{\mu_k}}
 >0.$$

 Next suppose that (\ref{step2-2}) holds.
 Substituting (\ref{lem5-1-1}) into (\ref{step2-3}) leads to
 \begin{eqnarray*}
 W'(t)&\le& \frac{\sum_{k \in J(t)}\frac{l_k}{\mu_k}(\rho_0+\frac{\lambda_k}{l_k}\frac{\zeta+\zeta^*}{2}
 +\lambda_2\frac{\zeta-\zeta^*}{2}-1)}
 {\frac{\zeta+\zeta^*}{2}+\sum_{k\in J(t)}\frac{l_k}{\mu_k}}.
 \end{eqnarray*}
 Hence $W'(t)\le -\epsilon$, where
 $$\epsilon \equiv \min_{K\in \{\{1\},\{3\},\{1,3\}\}}
 \frac{\sum_{k \in K}\frac{l_k}{\mu_k}(1-\rho_0-\frac{\lambda_k}{l_k}\frac{\zeta+\zeta^*}{2}
 -\lambda_2\frac{\zeta-\zeta^*}{2})}
 {\frac{\zeta+\zeta^*}{2}+\sum_{k\in K}\frac{l_k}{\mu_k}}
 >0.$$
 \qed

 \vspace{0.2cm}
 \noindent{\it Step 3.  There is a $\delta>0$ such that,
 for each fluid model solution $ \bar{\mathbb{X}}$ with
 $|\bar{Q}(0)| \le 1$, $W(t)=0$ for $t \ge \delta$.}

 \vspace{0.2cm}
 \noindent{\it Proof. }
 Suppose that $ \bar{\mathbb{X}}$ is a fluid model solution with
 $|\bar{Q}(0)| \le 1$.
 Then $W(0)\le \max_{1\le k \le 3}\frac{1}{\mu_k}$.
 According to (\ref{thm4-1}) and  (\ref{thm4-3}),
 we have
  $$W(\delta_m)\le \max_{1\le k \le 3}\frac{1}{\mu_k}+\rho_0
  \delta_m.$$
 By Steps 1 and 2,
  $$ W'(t)\le -\epsilon, ~~ \mbox{ if } t \ge \delta_m \mbox{ and }
  W(t)>0.$$
  Hence $W(t)=0$ for $t\ge \delta$, where $\delta\equiv
  \delta_m+\frac{1}{\epsilon}(\max_{1\le k \le 3}\frac{1}{\mu_k}+\rho_0
  \delta_m)$.
 \qed

 \vspace{0.4cm}
 \noindent {\it Proof of Proposition \ref{prop2}.}
 ~ Suppose that at least one of (\ref{inst1-1}), (\ref{inst2-1})
 or (\ref{inst3-1}) holds.
 Let $\bar {\mathbb{X}}$ be an augmented fluid model solution with
 $\bar Q(0)=0$.
 By (\ref{thm4-1}), (\ref{lem6-1}) and (\ref{lem6-2}), we have
 \begin{eqnarray*}
 \bar Q_2'(t)
 &>& \lambda_2-l_2 \frac{1-\rho_0+\beta_2}{\zeta+\frac{l_2}{\mu_2}} \nonumber \\
 &=&
 \frac{l_2(\frac{\lambda_2}{l_2}\zeta-1+\rho_0)}{\zeta+\frac{l_2}{\mu_2}},
 \label{prop2-pf1}
 \\
  \bar Q_j'(t)
 &>& \lambda_j-l_j \frac{1-\rho_0+\beta_j-\frac{\lambda_2}{l_2}\frac{\zeta-\zeta^*}{2}}{\frac{\zeta+\zeta^*}{2}+\frac{l_j}{\mu_j}} \nonumber \\
 &=&
 \frac{l_j(\frac{\lambda_j}{l_j}\frac{\zeta+\zeta^*}{2}+\frac{\lambda_2}{l_2}\frac{\zeta-\zeta^*}{2}-1+\rho_0)}
 {\frac{\zeta+\zeta^*}{2}+\frac{l_j}{\mu_j}}, ~~j=1,3.
 \label{prop2-pf2}
 \end{eqnarray*}
 Hence there exists $j\in \{1,2,3\}$ such that
 $\bar Q_j'(t)>0$ for all $t \ge 0$, which implies that
 the augmented fluid model is weakly unstable.
 \qed

 \vspace{0.4cm}
 \noindent {\it Proof of Proposition \ref{prop3}.}
 ~ Suppose that at least one of the
inequalities (\ref{inst1-1}), (\ref{inst2-1})
 or (\ref{inst3-1}) is strict.
Let $\bar {\mathbb{X}}$ be a fluid model solution with
 $\bar Q(0)=0$.

 First suppose that $\frac{\lambda_2}{l_2}>\frac{\lambda_k}{l_k}$,
 $k=1,3$.
 By (\ref{thm4-1}) and (\ref{lem5-4}), we have
 \begin{eqnarray*}
 \bar Q_2'(t)
 &\ge& \lambda_2-l_2 \frac{1-\rho_0+\sum_{k\in J(t)\cup \{2\}}\beta_k}{\zeta+\sum_{k\in J(t)\cup \{2\}}\frac{l_k}{\mu_k}} \nonumber \\
 &\ge&
 \frac{l_2(\frac{\lambda_2}{l_2}\zeta-1+\rho_0)}{\zeta+\sum_{k\in J(t)\cup
 \{2\}}\frac{l_k}{\mu_k}} \\
 &>&0,
 \end{eqnarray*}
 which implies that
 the fluid model is weakly unstable.

 Next suppose that there exists $j\in \{1,3\}$ such that
 $\frac{\lambda_j}{l_j}\ge\frac{\lambda_k}{l_k}$,
 $k=1,2,3$.
 By (\ref{thm4-1}) and (\ref{lem5-5}), we have
 \begin{eqnarray*}
  \bar Q_j'(t)
 &\ge& \lambda_j-l_j \dfrac{1-\rho_0+\sum_{k\in J(t)\cup\{j\}}\beta_k-\frac{\lambda_2}{l_2}\frac{\zeta-\zeta^*}{2}}
 {\frac{\zeta+\zeta^*}{2}+\sum_{k\in J(t)\cup\{j\}}\frac{l_k}{\mu_k}} \nonumber \\
 &\ge&
 \frac{l_j(\frac{\lambda_j}{l_j}\frac{\zeta+\zeta^*}{2}+\frac{\lambda_2}{l_2}\frac{\zeta-\zeta^*}{2}-1+\rho_0)}
 {\frac{\zeta+\zeta^*}{2}+\sum_{k\in J(t)\cup\{j\}}\frac{l_k}{\mu_k}} \\
 &>&0,
 \end{eqnarray*}
 which implies that
 the fluid model is weakly unstable.
 \qed

\section{Discussion on the fluid model for limited service disciplines
and simulation results} \label{section5}

 The upper and lower bounds for the fluid model equation
 (\ref{thm4-10}) is not tight if $l_2>1$.
 From the fluid model equations (\ref{thm4-1})-(\ref{thm4-13}), it
 can be observed that for each fluid model solution $\bar
 {\mathbb{X}}$, $\bar Q'(t)$ is determined by  $\bar U'_4(t)$,
 $J(t)$,
 and the model parameters $\lambda_j$, $\mu_j$, $l_j$, $j=1,2,3$,
 and $\nu_k$, $k=1,2,3,4$.

 For the case $l_2>1$, we conjecture that $\bar{U}^\prime_4(t)$ is
 not determined by only $J(t)$ and the model parameters but may also
depend on the
 distributions of the driving sequences, the
 inter-arrival, the service and the switch-over times.

To justify our conjecture, we present here two sets of simulation
results for the fluid limits in the system with limited service
disciplines. Specifically we observe fluid limits with $\bar
Q_2(t)=0$, $\bar Q_1(t)>0$ and $\bar Q_3(t)=0$.

 We remind that the condition $\bar Q_1(t)>0$ in a fluid limit
 corresponds to the condition in the real system that
  the first queue is infinitely large and that during each visit
to this queue the server serves exactly $l_1$ customers. The goal
is to find, in the long run, the fraction of time, $u_4$, the
server is in the switch-over regime from queue 1 to queue 3, and
the fraction of the number of the reduced cycles, $p$, among all.

\subsection{First set of examples}

We present five simulation Examples where we vary only one of the
distributions of the system (i.e. that of the inter-arrival times
to queue 2).

In the first three Examples, we are keeping the first moment fixed
and show that both fractions may differ significantly. These
examples illustrate that the stability conditions in the system
under consideration are not determined, in general, by the first
moments of the distributions of the primitive sequences.

We go further and show by example that the knowledge of the first
two moments is also not enough. We present Example 4 where the
distribution of inter-arrival times to queue 2 have the same 1st and
2nd moments with exponential distribution from Example 3, but here
also the fractions of interest significantly differ.

Finally, we finish with showing that even the knowledge of the
first three moments is insufficient. In Example 5, the
distribution of inter-arrival times to queue 2 have the same the
first, the second and the third moments with exponential
distribution from Example 3, but again with different fractions of
interest.

So, our conjecture is that the fractions of interest may depend on
the entirely whole distributions of the driving sequences, and
the knowledge of any finite number of moments is not sufficient to
determine the stability region precisely.

We consider the following system parameters:
$$
\lambda_2=\lambda_3=\frac14, \ \nu_1=2, \ \nu_3=1, \ \nu_4=3, \
\nu_2=+\infty, \ \mu_2=1, \ \mu_3=\frac23, \  l_2=4, \ l_3=2.
$$
We let for simplicity $\nu_2=+\infty$ that means that all
switch-over times $\xi_2(n)$ are zeros. We further assume for
simplicity that all $\sigma_1(n)=0$ too.

We recall that $u_4$ is the fraction of time when the server is
switching from queue 1 to queue 3, and that $p$ is the fraction of
the reduced cycles. In each of the following example, we run more
than $10^8$ cycles and find $u_4$ and $p$ with the error smaller
than $2\cdot 10^{-4}$ with probability greater than $0.9999$.

In what follows, a random variable with parameter, say $C$,
``has a uniform distribution'' means that it ``has a uniform
 distribution on the interval $(0,\, 2/C)$''.

 \begin{example}
 We assume that all interarrival, service, and switch-over
times are uniformly distributed. We obtain $u_4\approx 0.0466$ and
$p\approx 0.1825$.
 \end{example}

 \begin{example}
 We assume that the interarrival times to queue 2 have a
 probability density function $f(x)=8/x^3$, $x\ge 2$.
  The other primitive random variables are assumed to be uniformly
  distributed.
  Then we obtain
$u_4\approx 0,0518$, $p\approx 0,2027$.
 \end{example}

  \begin{example}
 We assume that the interarrival times to queue 2 have
 an exponential distribution with mean 4.
  The other primitive random variables are assumed to be uniformly
  distributed.
 Then we obtain $u_4\approx 0,0619$ and $p\approx 0,2410$.
 \end{example}

  \begin{example}
   We assume that the interarrival times to queue 2 have a
 probability density function $f(x)=ba^b/x^{b+1}$,
 $x\geq a$ with $a=8-4\sqrt{2}$ and $b=1+\sqrt{2}$.
  The other primitive random variables are assumed to be uniformly
  distributed.
 Then we obtain $u_4\approx 0,0446$ and $p\approx 0,1751$.
 \end{example}

 We remark that the first two moments (4 and 32 respectively) of the interarrival times to queue 2
 coincide in Examples 3 and 4.
 However  $u_4$ and $p$ are significantly different.

  \begin{example}
We assume that the interarrival times for queue 2 have a discrete
distribution given by
$$
\mathbf P( \tau_2(1) = 4(2-\sqrt{2}) ) = \frac{2+\sqrt{2}}{4},
\quad \mathbf P( \tau_2(1) = 4(2+\sqrt{2}) ) =
\frac{2-\sqrt{2}}{4}.
$$
  The other primitive random variables are assumed to be uniformly
  distributed.
 Then we obtain $u_4\approx 0.0641$ and $p\approx 0.2494$.
 \end{example}

 We remark that the first three moments (4, 32 and 384 respectively) of the interarrival times to queue 2
 coincide in Examples 3 and 5.
 However $u_4$ and $p$ are different.

\subsection{Examples with Weibull distribution}

In this subsection, we again vary only the distribution of the
interarrival times to queue 2. We focus on the class of Weibull
distributions with a fixed mean. More precisely, we assume that
the tail distribution function of the interarrival times for queue
2 is given by ${\mathbf P} (\tau_2(1)>x) = \exp (-bx^{a})$, $x\ge
0$, with $b=(\Gamma (1+a^{-1})/4)^a$. The other primitive random
variables are assumed to be uniformly
  distributed. Note that ${\mathbf E}
\tau_2(1)=4$. Let $l_2=6$ and $l_3=4$.

\begin{figure}
\includegraphics[scale=1]{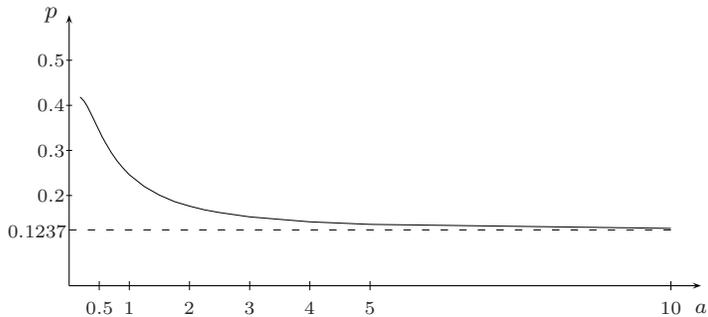}
\caption{Simulation results for $p$,
 varying the parameter $a$ for the Weibull distribution.} \label{fig1}
\end{figure}

 Table \ref{table1}  and Fig. \ref{fig1} present the simulation results for $p$ and $u_4$,
 varying the parameter $a$ for the Weibull distribution.
 The limiting value $0.1237$ in Fig. \ref{fig1} corresponds to
 the value of $p$ when the interarrival times for queue 2
 are 4, i.e., $\tau_2(n)=4$, $n=1,2,\ldots$, with probability 1.
 We observe that the lighter is the tail of the Weibull distribution,
 the smaller are the values $p$ and $u_4$.

 \begin{table}
 \begin{center}
 \begin{tabular}{c|c|crc|c|crc|c|c}
 \cline{1-3}\cline{5-7}\cline{9-11}
 $a$ &   $p$  & $u_4$  &\hspace{3mm}&
 $a$ &   $p$  & $u_4$  &\hspace{3mm}&
 $a$ &   $p$  & $u_4$
 \cr\cline{1-3}\cline{5-7}\cline{9-11}
 0.18 & 0.4181 & 0.1097 && 0.50 & 0.3435 & 0.0892 && 1.5  & 0.2009 & 0.0514 \\
 0.19 & 0.4174 & 0.1095 && 0.55 & 0.3297 & 0.0855 && 2    & 0.1765 & 0.0450 \\
 0.20 & 0.4162 & 0.1091 && 0.6 & 0.3179 & 0.0825 && 2.5  & 0.1623 &  0.0413 \\
 0.25 & 0.4089 & 0.1071 && 0.7 & 0.2953 & 0.0763 && 3    & 0.1527 &  0.0388 \\
 0.30 & 0.3982 & 0.1041 && 0.8 & 0.2762 & 0.0712 && 4    & 0.1417 & 0.0360 \cr
 0.35 & 0.3849 & 0.1005 && 0.9 & 0.2602 & 0.0670 && 5    & 0.1361 & 0.0346 \cr
 0.40 & 0.3713 & 0.0969 && 1 & 0.2461 & 0.0632 && 10   & 0.1272 & 0.0322 \cr
 0.45 & 0.3571 & 0.0929 && 1.25 & 0.2198 & 0.0563 && 20   & 0.1245 &  0.0315 \cr
 \cline{1-3}\cline{5-7}\cline{9-11}
\end{tabular}
\end{center}
 \caption{Simulation results for $p$ and $u_4$,
 varying the parameter $a$ for the Weibull distribution.
 }\label{table1}
 \end{table}

 \section{The fluid model and proof of Theorems \ref{th1} and \ref{th1-1}
for gated and exhaustive service disciplines} \label{section6}

 In this section, we consider the polling system with either gated
 or exhaustive service disciplines.
 The fluid limits and the fluid limits with fixed initial states
 are defined as in Section \ref{section3}.
 The following lemma can be shown by the standard procedure.

 \begin{lemma}\label{th2-1}
 For every fluid limit
 $\bar {\mathbb{X}}(t)=(\bar Q(t), \bar  T(t), \bar U(t))$,
 the following equations are satisfied:
 \begin{eqnarray}
 && \bar Q_k(t)=\bar Q_k(0)+\lambda_k t - \mu_k \bar T_k(t), ~~
    k=1,2,3, ~ t \ge 0; \label{thm4-1-1} \\
 && \bar Q_k(t)\ge 0, ~~ k=1,2,3, ~ t \ge 0; \label{thm4-2-1}\\
 && \bar T_k(\cdot) \mbox{ and } \bar U_j(\cdot) \mbox{ are nondecreasing},
    ~~ k=1,2,3,~j=1,2,3,4;  \label{thm4-3-1}\\
 && \sum_{k=1}^3\bar T_k(t)+\sum_{j=1}^4\bar U_j(t)=t, ~~~ t \ge 0; \label{thm4-4-1} \\
 && \mbox{If } |\bar Q(t)|>0, \mbox{ then } \sum_{k=1}^3\bar
 T_k'(t)=1, ~~ t\ge 0. \label{thm4-5-1}
 \end{eqnarray}
 \end{lemma}

 For the polling system with either gated
 or exhaustive service disciplines,
 we call the equations (\ref{thm4-1-1})-(\ref{thm4-5-1}) the {\it fluid model
 equations} and call a solution $ \bar{\mathbb{X}} =
 \{(\bar{Q}(t),\bar{T}(t),\bar{U}(t)): t \ge 0\}$, of the
 fluid model equations a {\it fluid model solution}.

 Using the similar argument as in the proof of Lemma \ref{th2-aug}, we
 can prove the following result.
  \begin{lemma}\label{th2-aug-1}
  Suppose that the Markov process $X$ is positive
 Harris
 recurrent, and
 a fluid limit with fixed initial state $\bar {\mathbb{X}}(t)=(\bar Q(t), \bar  T(t), \bar U(t))$
 is driven from a sample path $\omega$ that satisfies (\ref{f_i}).
 Then we have
 \begin{eqnarray}
 {\bar U}^\prime_k(t) &>&0, ~~t\ge 0, ~ k=1,2,3,4.\label{aug1-1-1}
 \end{eqnarray}
 \end{lemma}

   We call the equations (\ref{thm4-1-1})-(\ref{thm4-5-1}) plus (\ref{aug1-1-1})
  the {\it augmented fluid model
 equations} and call a solution $ \bar{\mathbb{X}}$, to these equations
 an {\it augmented fluid model solution}.

 Definitions \ref{def1} and \ref{def2}, and Lemmas \ref{lem-a1},
 \ref{lem-a2} and \ref{lem-a3} can be applied to the polling
 systems with gated and exhaustive service disciplines.
 Therefore Theorems \ref{th1} and \ref{th1-1} 
 for gated and exhaustive service disciplines
 are proved by Propositions \ref{prop4}, \ref{prop5} and
 \ref{prop6} below.
 \begin{proposition} \label{prop4}
 For the polling system with gated or exhaustive service
 disciplines, the fluid model is stable if $\rho_0<1$.
 \end{proposition}
 \proof
  For a fluid model solution
 $\bar{\mathbb{X}}$, let
 $W(t) \equiv \sum^3_{k=1}\frac{\bar Q_k(t)}{\mu_k}$.
 By (\ref{thm4-1-1}) and (\ref{thm4-5-1}), $W'(t)=\rho_0-1$ if
 $W(t)>0$. Hence the fluid model is stable if $\rho_0<1$.
 \qed

 \begin{proposition} \label{prop5}
 For the polling system with gated or exhaustive service
 disciplines, the augmented fluid model is weakly unstable if $\rho_0\ge 1$.
 \end{proposition}
 \proof
  For an augmented  fluid model solution
 $\bar{\mathbb{X}}$, let
 $W(t) \equiv \sum^3_{k=1}\frac{\bar Q_k(t)}{\mu_k}$.
 By (\ref{thm4-1-1}), (\ref{thm4-4-1}) and (\ref{aug1-1-1}), $W'(t)>\rho_0-1$.
 Hence the augmented fluid model is weakly unstable if $\rho_0\ge 1$.
 \qed

 \begin{proposition} \label{prop6}
 For the polling system with  gated or exhaustive service
 disciplines, the fluid model is weakly unstable if $\rho>1$.
 \end{proposition}
 \proof
  For a fluid model solution
 $\bar{\mathbb{X}}$, let
 $W(t) \equiv \sum^3_{k=1}\frac{\bar Q_k(t)}{\mu_k}$.
 By (\ref{thm4-1-1}), (\ref{thm4-3-1}) and (\ref{thm4-4-1}), $W'(t)\ge \rho_0-1$.
 Hence the fluid model is weakly unstable if $\rho>1$.
 \qed

\vskip2ex

The research of N.\,Chernova was partially supported by the
Ministry of Higher Education and Science of the Russian Federation
Grant~RNP.2.1.1.346. The research of S.\,Foss was partially
supported by the London Mathematical Society travel grant and
by the Royal Society International Joint Project. The research of
B. Kim was supported by the Korea Research Foundation (KRF) grant
funded by the Korea government
 (MEST) (2009-0076600).

\end{document}